\documentclass[12pt,oneside,reqno]{scrartcl}
\usepackage{cihan}
\author{C\.{i}han Bahran}
\date{}
\usepackage{etoolbox}
 \usepackage{relsize}
\usepackage[margin=1.0in]{geometry}
 
 \usepackage{mathtools}

\DeclarePairedDelimiter\floor{\lfloor}{\rfloor}
 
\newtheorem{thm-short}{Theorem}
\newtheorem{lem-short}[thm-short]{Lemma}
\newtheorem{cor-short}[thm-short]{Corollary}
\newtheorem{prop-short}[thm-short]{Proposition}

\title{\vspace{-1.0in}An improvement in the linear stable ranges for ordered configuration spaces}

\DeclareMathOperator{\conf}{PConf}
\DeclareMathOperator{\FI}{\mathbf{FI}}

\DeclareMathOperator{\co}{H}

\newcommand{\weak}{\delta}

\newcommand{\local}{h^{\text{max}}}

\newcommand{\confix}[2]{
  \ifstrempty{#1}{
    \conf_{\bullet}(#2)
  }
  {
    \conf_{#1}(#2)
  }
}

\newcommand{\man}{\mathcal{M}}

\AtBeginDocument{%
   \def\MR#1{}
}

\begin{document}
\maketitle
\vspace{-1.3cm}
In a recent paper, Church--Miller--Nagpal--Reinhold \cite[Application A, Theorem 4.3]{cmnr-range} gave the first linear stable ranges (in the sense of representation stability) for the \textbf{integral} cohomology of \textbf{ordered} configuration spaces of manifolds. When the manifold is orientable and is of dimension at least 3, the relevant spectral sequence of Totaro \cite{totaro-config} becomes sparse. We note that this sparsity can be exploited, while using no more than the methods in \cite{cmnr-range}, to improve the constants in the linear stable ranges. This note should be read after \cite{cmnr-range}, for the notation is borrowed from there.

\begin{thm-short} \label{main}
Let $\man$ be an \textbf{orientable}, connected manifold of dimension $d \geq 3$. Let $A$ be an abelian group. Then we have: 
\begin{birki}
 \item The stable degree $\weak(\co^{k}(\confix{}{\man};A))$ is $\leq k$.
 \vspace{1.5mm}
 \item The local degree \ds{\local(\co^{k}(\confix{}{\man};A))} is \ds{\leq \max \left\{-1, \, \frac{2d}{d-1}k  - 4\right\}}.
  \vspace{1.5mm}
 \item The generation degree \ds{t_{0}(\co^{k}(\confix{}{\man};A))} is \ds{\leq \max \left\{k, \, \frac{3d-1}{d-1}k - 3 \right\}}.
  \vspace{1.5mm}
 \item The presentation degree \ds{t_{1}(\co^{k}(\confix{}{\man};A))} is \ds{\leq \max \left\{k, \, \frac{5d-1}{d-1}k - 6\right\}}.
\end{birki}
\end{thm-short}
For the class of manifolds in the hypothesis of Theorem \ref{main}, instead of the coefficients $\frac{2d}{d-1}, \frac{3d-1}{d-1}, \frac{5d-1}{d-1}$ that we have in the items (2),(3),(4), Church--Miller--Nagpal--Reinhold \cite[Theorem 4.3]{cmnr-range} has the coefficients $4,5,9$, respectively. The item (1) is  included with no change for completeness.

A sequence of non-negatively graded cochain complexes $\{V^{\bullet}_{r}\}_{r=1}^{\infty}$ in an abelian category is called a (cohomological) \textbf{single-graded spectral sequence} if $\co^{k}(V^{\bullet}_{r}) \cong V^{k}_{r+1}$ for all $k \geq 0$ and $r \geq 1$. If the differentials 
$
 V_{r}^{k-1} \rarr V_{r}^{k} \rarr V_{r}^{k+1}
$
into and out of $V_{r}^{k}$ are zero when $r \geq r_{0}$, this implies 
$
 V_{r_{0}}^{k} = V_{r_{0} + 1}^{k} = \cdots 
$, 
and we write $V_{\infty}^{k}$ for this common value. We say the single-graded spectral sequence $\{V_{r}^{\bullet}\}$ \textbf{converges to} the graded object $M^{*}$ if $M^{*}$ has a filtration whose associated graded is $V^{\bullet}_{\infty}$.

The following lemma keeps consistent indexing with the proof of Church--Ellenberg--Farb \cite[Lemma 6.3.2]{cef}.
\begin{lem-short} \label{make-single}
 Let $\{E_{*}^{p,*}\}$ be a cohomological first quadrant spectral sequence bigraded in the standard way which converges to $M^{*}$. If there exists $D \geq 2$ such that $E_{2}^{p,*} = 0$ unless $* = qD$, then setting
\begin{align*}
  V^{k}_{r} := \bigoplus_{p + qD = k} E^{p,qD}_{rD + 1}
\end{align*} 
defines a single-graded spectral sequence $\{V_{r}^{\bullet}\}_{r=1}^{\infty}$ which converges to $M^{*}$. In addition, $V_{\infty}^{k} = V_{\floor*{k/D} + 1}^{k}$ for every $k$.
\end{lem-short}
\begin{proof}
 The vanishing assumption yields that the only nontrivial differentials occur on $(rD+1)$-th pages, with $E_{(r-1)D + 2} = \cdots E_{rD} = E_{rD+1}$. Thus we do have a single-graded spectral sequence $\{V_{r}^{\bullet}\}$ as described.
 
 For the last claim, first note that the differentials on the $(rD+1)$-th page of $\{E_{*}^{p,*}\}$ are of the form 
\begin{align*}
 E_{rD+1}^{p-rD-1, (q+r)D} \rarr E_{rD+1}^{p,qD} \rarr E_{rD+1}^{p+rD+1, (q-r)D} \, .
\end{align*}
Observe that whenever $p+qD = k$ (with $p,q \geq 0$) and $r \geq \floor*{k/D} + 1$,
\begin{itemize}
 \item $p-rD-1 < p - k -1 < 0$,
 \item $(q-r)D = qD - rD < qD - k \leq 0$.
\end{itemize}
Therefore $E_{rD+1}^{p,qD} = E_{\infty}^{p,qD}$ whenever $p+qD = k$ and $r \geq \floor*{k/D} + 1$. 
\end{proof}

\begin{prop-short} [{\cite[Proposition 4.1]{cmnr-range}}] \label{range}
 Let $\{V^{\bullet}_{r}\}$ be a single-graded spectral sequence of $\FI$-modules presented in finite degrees, such that for each $k$ the $\FI$-module $V^{k}_{1}$ is semi-induced and generated in degree $D_{k}$. Then we have
\begin{enumerate}
 \item $\weak(V_{r}^{k}) \leq D_{k}$,
 \item $\local(V_{r}^{k}) \leq \max\left( \left\{ 2 D_{\ell} - 2 : \ell \leq k + r - 2\right\} \cup \{-1\} \right)$,
\end{enumerate}
for every $r$.
\end{prop-short}

\begin{cor-short} \label{dur-burda}
 Let $\{E_{*}^{p,*}\}$ be a cohomological first quadrant spectral sequence of $\FI$-modules converging to $M^{*}$, such that $E_{2}^{p,*} = 0$ unless $* = qD$ for some $D \geq 2$. Suppose in addition that for all $p,q$, the $\FI$-module $E_{2}^{p,qD} = E_{D+1}^{p,qD}$ is semi-induced and generated in degree at most $D_{k}$ whenever $p+qD = k$. Then we have 
\begin{enumerate}
 \item $\weak(M^{k}) \leq D_{k}$,
 \item $\local(M^{k}) \leq \max \left\{ 2 D_{\ell} - 2 : \ell \leq k + k/D - 1\right\}$.
\end{enumerate}
\end{cor-short}
\begin{proof}
 By the first part of Lemma \ref{make-single}, setting 
\begin{align*}
  V^{k}_{r} := \bigoplus_{p + qD = k} E^{p,qD}_{rD + 1}
\end{align*}
yields single-graded spectral sequence $\{V^{\bullet}_{r}\}_{r=1}^{\infty}$ of $\FI$-modules that converges to $M^{*}$. It follows from the hypotheses that $V_{1}^{k}$ is semi-induced and generated in degree at most $D_{k}$. Thus, by the second part of Lemma \ref{make-single} and Proposition \ref{range}, we get 
\begin{enumerate}
 \item $\weak(V_{\infty}^{k}) \leq D_{k}$,
 \item $\local(V_{\infty}^{k}) \leq \max \left( \left\{ 2 D_{\ell} - 2 : \ell \leq k + k/D - 1\right\} \cup \{-1\} \right)$.
\end{enumerate}
The analogous claims for $M^{k}$ follow from \cite[Proposition 3.2]{cmnr-range}.
\end{proof}


\begin{proof}[Proof of Theorem \ref{main}]
 By Totaro's and Church's \cite{church-config} work (see \cite[Section 6.3]{cef} for details), there is a first quadrant spectral sequence $E_{*}^{p,*}$ of $\FI$-modules converging to $\co^{*}(\confix{}{\man};A))$ that satisfies the hypotheses of Corollary \ref{dur-burda} with $D = d-1$ and $D_{k} = k$. Thus (1) and (2) follow from Corollary \ref{dur-burda}. And then (3) and (4) follow from \cite[Proposition 3.1]{cmnr-range}.
\end{proof}

\bibliographystyle{hamsalpha}
\bibliography{stable-boy}

\end{document}